\documentclass[a4paper,12pt]{article}
\usepackage{pdfsync}
\usepackage{amsfonts}
\usepackage{amssymb}
\usepackage{amsmath}
\usepackage{mathptmx}
\usepackage{cancel}
\usepackage{pstricks}
\usepackage{psfrag}
\usepackage{mathrsfs}
\usepackage{graphicx}
\usepackage{pgf,tikz}
\usetikzlibrary{arrows}
\usepackage{algorithm}
\usepackage{algpseudocode}
\usepackage{soul}
\def\be{\begin{equation}}
\def\ee{\end{equation}}
\def\bea{\begin{eqnarray}}
\def\eea{\end{eqnarray}}
\def\beann{\begin{eqnarray*}}
\def\eeann{\end{eqnarray*}}

\def\ns{\hspace{-1mm}}

\def\spacingset#1{\def\baselinestretch{#1}\small\normalsize}
\spacingset{1.30}

\newtheorem{lemma}{Lemma}
\newtheorem{theorem}{Theorem}
\newtheorem{remark}{Remark}

\newtheorem{proposition}{Proposition}
\newtheorem{definition}{Definition}

\def\be{\begin{equation}}
\def\ee{\end{equation}}
\def\bea{\begin{eqnarray}}
\def\eea{\end{eqnarray}}
\def\beann{\begin{eqnarray*}}
\def\eeann{\end{eqnarray*}}

\def\ns{\hspace{-1mm}}

\def\proof{\noindent{\bf{\em Proof:}\ \ }}
\def\QED{\mbox{\rule[0pt]{1.5ex}{1.5ex}}}
\def\endproof{\hspace*{\fill}~\QED\par\endtrivlist\unskip}
\def\t{\textrm}

\newcommand{\defi}{\stackrel{\text{\tiny def}}{=}}

\definecolor{Royalblue}{cmyk}{1,0.30,0.2,0.2}

\def\gC{{\cal C}}

\def\gE{{\cal E}}

\def\gP{{\cal P}}

\def\bmat{\left[ \begin{array}}
\def\emat{\end{array} \right]}

\def\bmat{\left[ \begin{array}}
\def\emat{\end{array} \right]}
\def\bsmat{\left[ \begin{smallmatrix}}
\def\esmat{\end{smallmatrix} \right]}

\def\gP{{\cal P}}

\definecolor{light-gray}{gray}{0.84}
\begin{document}
\begin{titlepage}
\title{\vspace{-5mm}
A new parameterisation of Pythagorean triples in terms of odd and even series\vspace{10mm}}
\author{Anthony Overmars, Lorenzo Ntogramatzidis}
\date{\small Department of Mathematics and Statistics,\\[-2pt]
         Curtin University, Perth (WA), Australia \\[-2pt]
         {\tt \{Anthony.Overmars,L.Ntogramatzidis\}@curtin.edu.au}         
                 }
\thispagestyle{empty} \maketitle \thispagestyle{empty}
\begin{abstract}%
{\normalsize
In this paper we introduce a formula that parameterises the Pythagorean triples as elements of two series. With respect to the standard Euclidean formula, this parameterisation does not generate the Pythagorean triples where the elements of the triple are all divisible by $2$. A necessary and sufficient condition is also proposed for a Pythagorean triple obtained from this formula to be primitive.
}
\end{abstract}

{\small
\begin{center}
\begin{minipage}{14.2cm}
\vspace{2mm}
\end{minipage}
\end{center}
}
\thispagestyle{empty}
\end{titlepage}
%

\section{Introduction}
\label{secintro} 
The Pythagorean theorem asserts that if $a$, $b$, $c$ are the sides of a right triangle, with $c$ being the hypotenuse, then $a^2 +b^2=c^2$. In the case where $a$, $b$ and $c$ are all natural numbers different from zero, the triple $(a, b, c)$ is called a {\em Pythagorean triple}. Therefore, finding Pythagorean triples is equivalent to finding right triangles with integral sides. The problem for rational numbers is clearly the same up to a scale factor. The notion of Pythagorean triple, and its relation to the Pythagorean theorem, is a cornerstone of several areas of pure mathematics, including number theory, elementary and algebraic geometry,  as well as applied mathematics, due to its relevance in areas such as cryptography and random number generation algorithms.
Given the large number of contributions in these areas, it would be impossible to quote even only a fraction of the relevant references, and we consequently direct the interested reader to the  monographs \cite{Maor-07} and \cite{Sierpinski}, which provide surveys of the extensive literature in these areas, and to the papers \cite{Cass-A-90,Frisch-V-08,Vasertein-SF-10} along with the references cited therein.

Euclid's formula is a fundamental instrument for generating Pythagorean triples given an arbitrary pair of non-zero natural numbers $u,v$. Euclid's formula states that 
all primitive Pythagorean triples $(a,b,c)$ in which $b$ is even are obtained from the formulae
\[
a=u^2-v^2, \qquad b=2\,u\,v, \qquad c=u^2+v^2,
\]
where $u>v$, and $u$ and $v$ being all pairs of relatively prime numbers of which one of them (whatsoever) is even and the other is odd.
Each primitive Pythagorean triple $(a,b,c)$ where $b$ is even is obtained in this way only once.

In \cite{Harriot-70}, Thomas Harriot was the first to suggest that Pythagorean triples
exist in series, but no general form in terms of series was expressed therein. The purpose of this paper is to show that Pythagorean triples  can indeed be characterised in terms of odd and even series, which comprise the triangles in which the difference between the hypotenuse and a side is odd and/or even.
 The significance of odd and even sides was recognised and
expressed in \cite[Lemma 1]{Cass-A-90}, but, to the best of the authors' knowledge, a complete formula parameterising the Pythagorean triples in this way has not yet been proposed.

\section{Preliminary definitions and results}
A generic Pythagorean triple is denoted as the triple
\bea
\label{pt}
(a,b,c)
\eea
where $a,b,c \in \mathbb{N}\setminus \{0\}$, where $c> \max \{a,b\}$, i.e., $a$ and $b$ denote the sides of a right triangle, and $c$ denotes the length of the hypothenuse. The set of all Pythagorean triples is denoted by $\gP$.
A Pythagorean triple $(a,b,c)$ is said to be {\em primitive} if $a$, $b$ and $c$ are relatively prime. We denote by $\gP_{\scriptscriptstyle 0}$ the set of all primitive Pythagorean triples.

We recall the so-called Euclid's formula for parameterising Pythagorean triples given an arbitrary pair of natural numbers $u,v \in \mathbb{N}\setminus \{0\}$ with $u>v$.

\begin{lemma}{\sc [Euclid's Formula]}\\
\label{Leuclide}
For any $u,v \in \mathbb{N}\setminus \{0\}$ with $u>v$, 
\bea
\label{euclide}
\bigl(u^2-v^2,2\,u\,v,u^2+v^2\bigr)
\eea
is a Pythagorean triple. 
\end{lemma}

The Pythagorean triples generated by (\ref{euclide}) will be referred to as the {\em Euclidean triples}. The set of such triples is denoted by $\gE$. 
Euclid's formula generates all primitive Pythagorean triples. Indeed, the following important result holds, see e.g. \cite{Maor-07}.

\begin{proposition}
\label{prop1}
There holds $\gP \supset \gE \supset \gP_{\scriptscriptstyle 0}$. Moreover, a Euclidean triple generated with (\ref{euclide}) is primitive if and only if $u$ and $v$ are relatively prime and $u-v$ is odd.
\end{proposition}

The inclusion of $\gE \supset \gP_{\scriptscriptstyle 0}$ in Proposition \ref{prop1} 
shows that Euclid's formula generates all primitive Pythagorean triples. However, such inclusion 
is strict, in the sense that (\ref{euclide}) generates also a subset of $\gP \setminus \gP_{\scriptscriptstyle 0}$. 
For example, if $u$ and $v$ are chosen to be odd, then $a=u^2-v^2$, $b=2\,u\,v$ and $c=u^2+v^2$ are even, which means that the triple $(a,b,c)$ is not primitive.
If $u=6$ and $v=3$, they are not relatively prime, and $u-v$ is odd. In this case the corresponding Euclidean triple is $(27,36,45)$, which is not primitive.
Moreover, since the inclusion $\gP \supset \gE$ in Proposition \ref{prop1} is also strict,  $\gE$ is a proper subset of the set $\gP$ of all Pythagorean triples. Indeed, while it is true that every primitive Pythagorean triple is an element of $\gE$, some non-primitive Pythagorean triples are not in $\gE$. For example, for the non-primitive Pythagorean triple $(9,12,15) \in \gP \setminus \gP_{\scriptscriptstyle 0}$, values $u,v \in \mathbb{N}\setminus \{0\}$ with $u>v$ such that $\bigl(u^2-v^2,2\,u\,v,u^2+v^2\bigr)=(9,12,15)$ cannot be found.


 If $(a,b,c)$ is a primitive Pythagorean triple, then exactly one of $a$, $b$ is odd, and $c$ is odd, see e.g. \cite{Maor-07}.
 The converse is not necessarily true, i.e., there are Pythagorean triples that are not primitive for which one between $a$ and $b$ is odd and $c$ is odd, e.g. $(27,36,45)$.
 
In this paper, we are interested in a subset of the set of Euclidean triples, herein denoted with the symbol $\gC$, which comprises the Euclidean triples for which exactly one between $a$ and $b$ is odd and $c$ is odd. It is easily established that exactly one between $a$ and $b$ is odd if and only if exactly one between $u$ and $v$ is odd, so that
\beann
\gC \ns&\ns = \ns&\ns \{(a,b,c)\in \gE \,|\, \text{$a-b$ and $c$ are odd}\}\\
\ns&\ns = \ns&\ns \left\{\bigl(u^2-v^2,2\,u\,v,u^2+v^2\bigr)\,|\,u,v \in \mathbb{N}\setminus \{0\}, \; u>v, \; \text{$u-v$ is odd}\right\}.
\eeann
 From the considerations following Proposition \ref{prop1} it turns out that
 \[
\gP \supset \gE \supset \gC \supset \gP_{\scriptscriptstyle 0}.
 \]

It is assumed that the order of the elements of a triple in $\gC$ is the one induced by the order imposed in $\gE$. Thus, $a$ is odd and $b$ is even for any $(a,b,c)\in \gC$. 
Therefore, for every $(a,b,c)\in \gC$, the difference $c-b$ is odd and the difference $c-a$ is even.

\section{Odd and even series of Pythagorean triples}

The following preliminary result is important, because it enables the parameterisation of the set of Pythagorean triples in a convenient form.

\begin{lemma}
\label{fundamental}
For any Pythagorean triple $(a,b,c) \in \gC$
\begin{itemize}
\item there exists $m \in  \mathbb{N}\setminus \{0\}$ such that $c-b=(2\,m-1)^2$;
\item there exists $n \in  \mathbb{N}\setminus \{0\}$ such that $c-a=2\,n^2$.
\end{itemize}
Conversely, given $m,n \in \mathbb{N}\setminus \{0\}$, there exists $(a,b,c) \in \gC$ such that $c-b=(2\,m-1)^2$ and $c-a=2\,n^2$.
\end{lemma}
\proof
Firstly, it is shown that given $(a,b,c) \in \gC$, the difference $c-b$ can be written as $(2\,m-1)^2$ for some $m \in  \mathbb{N}\setminus \{0\}$. From Lemma \ref{Leuclide} and the definition of $\gC$, there exist $u,v \in \mathbb{N}\setminus \{0\}$ with $u>v$ and $u-v$ odd such that $a=u^2-v^2$, $b=2\,u\,v$ and $c=u^2+v^2$. We can write
\beann
c-b \ns&\ns = \ns&\ns (u^2+v^2)-2\,u\,v = (u-v)^2.
 \eeann
 Thus, $c-b$ can be written as $(2\,m-1)^2$ for some $m \in  \mathbb{N}\setminus \{0\}$ if and only if $u-v$ is odd, which is true.
 
Secondly, given $(a,b,c) \in \gC$,  the difference $c-a$ can be shown to be equal to $2\,n^2$ for some $n \in  \mathbb{N}\setminus \{0\}$. Indeed, using again Lemma \ref{Leuclide}, it is found that
 \beann
c-a \ns&\ns = \ns&\ns (u^2+v^2)-(u^2+v^2)=2\,v^2.
 \eeann
Then, $n=v$ is a solution. 
To prove the converse, in view of Lemma \ref{Leuclide}, it suffices to prove that the equations
\beann
u^2+v^2-2\,u\,v  \ns&\ns = \ns&\ns (2\,m-1)^2, \\
 u^2+v^2-(u^2-v^2)  \ns&\ns = \ns&\ns 2\,n^2, 
 \eeann
 can be solved in $u,v \in \mathbb{N}\setminus \{0\}$ with $u>v$, and where $u-v$ is odd.
It is straightforward to verify that $u = n +2\,m-1$ and $v=n$ is a solution.
\endproof 

Loosely, in view of Lemma \ref{fundamental}, for every Pythagorean triple $(a,b,c)\in \gC$, the difference $c-b$ is constrained to be an odd integer in $\{1,9,25,49,....\}$. Similarly, the difference $c-a$  is constrained to be equal to one of the even numbers in $\{2,8,18,32,....\}$. This simple consideration enables the introduction of the following definition.
\begin{definition}
For every $m,n \in  \mathbb{N}\setminus \{0\}$, the {\em series of odd Pythagorean triples} of $\gC$ is defined as
\[
{\rm odd}(m) \defi \left\{ (a,b,c) \in \gC\,|\,c-b=(2\,m-1)^2\right\}
\]
and the {\em series of even Pythagorean triples} of $\gC$ is defined as
\[
{\rm even}(n) \defi \left\{ (a,b,c) \in \gC\,|\,c-a=2\,n^2\right\}.
\]
\end{definition}
For example, $(3,4,5),(5,12,13),(7,24,25) \in {\rm odd}(1)$ and  $(3,4,5),(15,8,17),(35,12,37)\in {\rm even}(1)$.

As a consequence of Lemma \ref{fundamental}, the following results hold.

\begin{theorem}
\label{cor}
The following facts hold true: 
\begin{enumerate}
\item There holds
\[
\gC=\bigcup_{m \in  \mathbb{N}\setminus \{0\}} {\rm odd}(m) =\bigcup_{n \in  \mathbb{N}\setminus \{0\}} {\rm even}(n).
\]
\item For any $m_1,m_2\in  \mathbb{N}\setminus \{0\}$ with $m_1 \neq m_2$, there holds ${\rm odd}(m_1)\cap {\rm odd}(m_2)=\emptyset$. Likewise, for any $n_1,n_2\in  \mathbb{N}\setminus \{0\}$ with $n_1 \neq n_2$, there holds ${\rm even}(n_1)\cap {\rm even}(n_2)=\emptyset$. 
\item For every $(a,b,c)\in \gC$ there exists $n,m\in  \mathbb{N}\setminus \{0\}$ such that $(a,b,c)\in{\rm odd}(m) \cap {\rm even}(n)$. Conversely, for every  $n,m\in  \mathbb{N}\setminus \{0\}$ the intersection ${\rm odd}(m) \cap {\rm even}(n) \in \gC$ contains a single Pythagorean triple.
\item Let $m,n \in  \mathbb{N}\setminus \{0\}$. Then, the triple $(a,b,c)\in \gC$ such that $c-b=(2\,m-1)^2$ and $c-a=2\,n^2$ can be written in a unique way as 
 the Diophantine equations in $m,n$
\bea
a \ns&\ns = \ns&\ns -2\,n+4\,n\,m+4\,m^2-4\,m+1, \label{dioph1} \\
b \ns&\ns = \ns&\ns 2\,n^2-2\,n+4\,n\,m, \label{dioph2} \\
c \ns&\ns = \ns&\ns 2\,n^2-2\,n+4\,n\,m +4\,m^2-4\,m+1.\label{dioph3}
\eea
\end{enumerate}
\end{theorem}
\proof
Claims 1-3 are a simple consequence of Lemma \ref{fundamental}.
 We only need to prove the last point. From the last part of the proof of Lemma \ref{fundamental} it is found that
 \beann
a \ns&\ns = \ns&\ns (n+2\,m-1)^2-n^2,\\
b\ns&\ns = \ns&\ns 2\,(n+2\,m-1)\,n,\\
c\ns&\ns = \ns&\ns (n+2\,m-1)^2+n^2,
\eeann
which immediately lead to (\ref{dioph1}-\ref{dioph3}). It is now shown that $(a,b,c)\in \gC$ is written in a unique way as a function of $n,m \in  \mathbb{N}\setminus \{0\}$. To this end, using (\ref{dioph1}) yields
\bea
\label{nfromabc}
n=\frac{a-(2\,m-1)^2}{2\,(2\,m-1)},
\eea
which, once substituted in (\ref{dioph2}), gives the biquadratic equation
$(2\,m-1)^4+2\,b\,(2\,m-1)-a^2=0$, 
which is easily seen to admit only one positive solution 
\[
m=\frac{1+\sqrt{c-b}}{2}.
\]
Replacing the latter into (\ref{nfromabc}) yields
\bea
\label{nfromabc1}
n=\frac{a+b-c}{2\,\sqrt{c-b}}.
\eea
This shows that given $(a,b,c)\in \gC$, the values of the natural numbers $m$ and $n$ such that (\ref{dioph1}-\ref{dioph3}) hold are uniquely determined.
\endproof

\begin{remark}
{\em
Notice that (\ref{dioph1}-\ref{dioph3}) can be written as
\bea
\label{matrixform}
 \bmat{c} a \\ b \\ c \emat=
\bmat{ccc} 0 & 1 & 1 \\ 1 & 1 & 0 \\ 1 & 1 & 1 \emat \bmat{c} 2\,n^2 \\ 2\,n\,(2\,m-1) \\ (2\,m-1)^2 \emat.
\eea
Consider the right triangle in Figure \ref{triangle}.
\begin{figure}[!h]
\begin{center}
\begin{tikzpicture}
\draw [thick,light-gray,fill] (0,0) --  (3,0) -- (0,4) -- (0,0);
\draw [<->,thin] (3+0.2,0+0.2) --  (0+0.2,4+0.2);
\draw [<->,thin] (0,0-0.3) --  (3,0-0.3);
\draw [<->,thin] (0-0.3,0) --  (0-0.3,4);
\draw [thick] (0,0) --  (3,0) -- (0,4) -- (0,0);
\draw [thin,dashed,shift={(3,0)}] plot[domain=0:6.28,variable=\t]({1*1*cos(\t r)+0*1*sin(\t r)},{0*1*cos(\t r)+1*1*sin(\t r)});
\draw [thin,dashed,shift={(3,0)}] plot[domain=1:3.0,variable=\t]({1*3*cos(\t r)+0*3*sin(\t r)},{0*3*cos(\t r)+1*3*sin(\t r)});

\draw [thin,dashed,shift={(0,4)}] plot[domain=-3.23:0.3,variable=\t]({1*2*cos(\t r)+0*2*sin(\t r)},{0*2*cos(\t r)+1*2*sin(\t r)});

\draw [thin,dashed,shift={(0,4)}] plot[domain=2.5*1.7:1.57*4.0,variable=\t]({1*4*cos(\t r)+0*4*sin(\t r)},{0*4*cos(\t r)+1*4*sin(\t r)});

\draw[color=black] (-.05,4.3) node {\small $A$};
\draw[color=black] (3.2,-0.3) node {\small $B$};
\draw[color=black] (1.4,-0.5) node {\small $a$};
\draw[color=black] (-0.3,-0.3) node {\small $C$};
\draw[color=black] (-.6,1.8) node {\small $b$};

\draw[color=black] (2.1,2.2) node {\small $c$};

\draw[color=black] (2.45,0.4) node {\small $d$};
\draw[color=black] (1.5,1.5) node {\small $f$};
\draw[color=black] (0.5,3.0) node {\small $e$};
%
%
\end{tikzpicture}
\end{center}
\vspace{-1cm}
\caption{Geometric interpretation of $m$ and $n$.}
\label{triangle}
\end{figure}
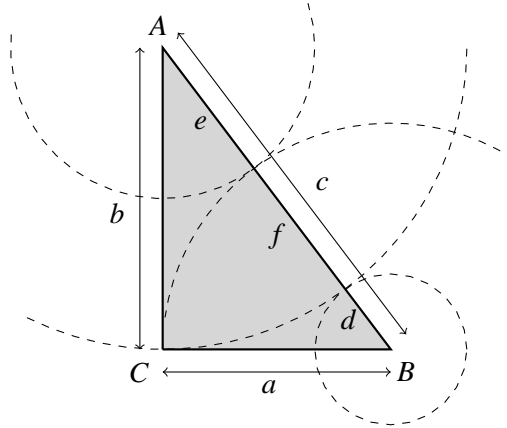
 Let $a=\overline{B\,C}$, $b=\overline{A\,C}$ and $c=\overline{A\,B}$.
Suppose $d=c-b$ is odd and $e=c-a$ is even. Then, defining $f=a-d$, the three relations
\beann
\begin{array}{rcl}
a \ns&\ns = \ns&\ns f+d \\ 
b\ns&\ns = \ns&\ns e+f \\
 c\ns&\ns = \ns&\ns e+d+f \end{array} \quad \Leftrightarrow \qquad 
 \bmat{c} a \\ b \\ c \emat=
\bmat{ccc} 0 & 1 & 1 \\ 1 & 1 & 0 \\ 1 & 1 & 1 \emat \bmat{c} e \\ f \\ d \emat
\eeann
are obtained.
From the non-singularity of $\bsmat 0 && 1 && 1 \\[1mm] 1 && 1 && 0 \\[1mm] 1 && 1 && 1 \esmat$ it follows that  $e=2\,n^2$ and $d=(2\,m-1)^2$. 
}
\end{remark}

\begin{remark}
\label{remChat}
{\em The Pythagorean triples obtained in this way are all such that $c$ is odd. However, 
 if $m$ is allowed to also take the values $1/2,3/2,5/2....$, then $a$, $b$, $c$ are still in $\mathbb{N} \setminus \{0\}$ and $(a,b,c)$ is still a Pythagorean triple, but $a$, $b$ and $c$ are all even. This means that $\gC$ comprises all the triples of $\gE$ except for those for which $a,b,c$ are even. 
 Stated differently, denoting $\mathbb{N}_{1/2} \defi \mathbb{N}\cup \left\{ p+\frac{1}{2} \,|\,p\in \mathbb{N} \right\}$, and defining 
 \[
 \widehat{\gC} \defi \bigcup_{m \in \mathbb{N}_{1/2}\setminus\{0\}} \text{odd}(m),
 \]
there holds $\widehat{\gC}=\gE$. 
 Indeed, defining $s=2\,m$, the equations in the proof of Theorem \ref{cor} become 
\beann
a \ns&\ns = \ns&\ns (n+s-1)^2-n^2, \qquad 
b= 2\,(n+s-1)\,n,\qquad
c= (n+s-1)^2+n^2,
\eeann
so that by setting $u=n+s-1$ and $v=n$ it is seen that $(a,b,c) \in \gE$.
}
\end{remark}

As a result of Lemma \ref{fundamental} and Theorem \ref{cor}, every Pythagorean triple $(a,b,c)$ in $\gC$ can be expressed in a unique way as a function $\mathfrak{P}$ of $m$ and $n$, i.e., we can write
\[
(a,b,c)=\mathfrak{P}(m,n).
\]
In view of the considerations above, all triples in $\gC$ can be represented in terms of $m$ and $n$ into a lattice as shown in Figure \ref{lattice1}.
\begin{center}
\begin{figure}[!h]
\hspace{1cm} \begin{tikzpicture}[line cap=round,line join=round,>=triangle 45,x=1.4cm,y=1.23cm]

\draw [->,thick] (-0.5,0) -- (6,0);
\draw [->,thick] (0,-0.5) -- (0,6);

\draw [fill] (0,0) circle [radius=0.02];
\draw [fill] (0,1) circle [radius=0.02];
\draw [fill] (0,2) circle [radius=0.02];
\draw [fill] (0,3) circle [radius=0.02];
\draw [fill] (0,4) circle [radius=0.02];
\draw [fill] (0,5) circle [radius=0.02];

\draw [fill] (1,0) circle [radius=0.02];
\draw [fill] (1,1) circle [radius=0.02];
\draw [fill] (1,2) circle [radius=0.02];
\draw [fill] (1,3) circle [radius=0.02];
\draw [fill] (1,4) circle [radius=0.02];
\draw [fill] (1,5) circle [radius=0.02];

\draw [fill] (2,0) circle [radius=0.02];
\draw [fill] (2,1) circle [radius=0.02];
\draw [fill] (2,2) circle [radius=0.02];
\draw [fill] (2,3) circle [radius=0.02];
\draw [fill] (2,4) circle [radius=0.02];
\draw [fill] (2,5) circle [radius=0.02];

\draw [fill] (3,0) circle [radius=0.02];
\draw [fill] (3,1) circle [radius=0.02];
\draw [fill] (3,2) circle [radius=0.02];
\draw [fill] (3,3) circle [radius=0.02];
\draw [fill] (3,4) circle [radius=0.02];
\draw [fill] (3,5) circle [radius=0.02];

\draw [fill] (4,0) circle [radius=0.02];
\draw [fill] (4,1) circle [radius=0.02];
\draw [fill] (4,2) circle [radius=0.02];
\draw [fill] (4,3) circle [radius=0.02];
\draw [fill] (4,4) circle [radius=0.02];
\draw [fill] (4,5) circle [radius=0.02];

\draw [fill] (5,0) circle [radius=0.02];
\draw [fill] (5,1) circle [radius=0.02];
\draw [fill] (5,2) circle [radius=0.02];
\draw [fill] (5,3) circle [radius=0.02];
\draw [fill] (5,4) circle [radius=0.02];
\draw [fill] (5,5) circle [radius=0.02];

%
\draw [dashed,thin] (-0.5,1) -- (6,1);
\draw [dashed,thin] (-0.5,2) -- (6,2);
\draw [dashed,thin] (-0.5,3) -- (6,3);
\draw [dashed,thin] (-0.5,4) -- (6,4);
\draw [dashed,thin] (-0.5,5) -- (6,5);

\draw [dashed,thin] (1,-0.5) -- (1,6);
\draw [dashed,thin] (2,-0.5) -- (2,6);
\draw [dashed,thin] (3,-0.5) -- (3,6);
\draw [dashed,thin] (4,-0.5) -- (4,6);
\draw [dashed,thin] (5,-0.5) -- (5,6);
\node [above right] at (1,1) {\tiny $\!\!\!\!(3,4,5)$};
\node [above right] at (2,1) {\tiny $\!\!\!\!(15,8,17)$};
\node [above right] at (3,1) {\tiny $\!\!\!\!(35,12,37)$};
\node [above right] at (4,1) {\tiny $\!\!\!\!(63,16,65)$};
\node [above right] at (5,1) {\tiny $\!\!\!\!(99,20,101)$};

\node [above right] at (1,2) {\tiny $\!\!\!\!(5,12,13)$};
\node [above right] at (2,2) {\tiny $\!\!\!\!(21,20,29)$};
\node [above right] at (3,2) {\tiny $\!\!\!\!(45,28,53)$};
\node [above right] at (4,2) {\tiny $\!\!\!\!(77,36,85)$};
\node [above right] at (5,2) {\tiny $\!\!\!\!(117,44,125)$};

\node [above right] at (1,3) {\tiny $\!\!\!\!(7,24,25)$};
\node [above right] at (2,3) {\tiny $\!\!\!\!(27,36,45)$};
\node [above right] at (3,3) {\tiny $\!\!\!\!(55,48,73)$};
\node [above right] at (4,3) {\tiny $\!\!\!\!(91,60,109)$};
\node [above right] at (5,3) {\tiny $\!\!\!\!(135,72,153)$};

\node [above right] at (1,4) {\tiny $\!\!\!\!(9,40,41)$};
\node [above right] at (2,4) {\tiny $\!\!\!\!(33,56,65)$};
\node [above right] at (3,4) {\tiny $\!\!\!\!(65,72,97)$};
\node [above right] at (4,4) {\tiny $\!\!\!\!(105,88,137)$};
\node [above right] at (5,4) {\tiny $\!\!\!\!(153,104,185)$};

\node [above right] at (1,5) {\tiny $\!\!\!\!(11,60,61)$};
\node [above right] at (2,5) {\tiny $\!\!\!\!(39,80,89)$};
\node [above right] at (3,5) {\tiny $\!\!\!\!(75,100,125)$};
\node [above right] at (4,5) {\tiny $\!\!\!\!(119,120,169)$};
\node [above right] at (5,5) {\tiny $\!\!\!\!(171,140,221)$};

\node [below] at (6,0) { $m$};
\node [left] at (0,6) { $n$};
\draw [-,thick] (0,-1) -- (3.5,6);
\node [right] at (-0.4,-1.2) {\small $n=2\,m-1$};

\node [below right] at (1,0) {\tiny $1$};
\node [below right] at (2,0) {\tiny $2$};
\node [below right] at (3,0) {\tiny $3$};
\node [below right] at (4,0) {\tiny $4$};
\node [below right] at (5,0) {\tiny $5$};

\node [below left] at (0,1) {\tiny $1$};
\node [below left] at (0,2) {\tiny $2$};
\node [below left] at (0,3) {\tiny $3$};
\node [below left] at (0,4) {\tiny $4$};
\node [below left] at (0,5) {\tiny $5$};

\node [below] at (1,-0.5) {\small odd($1$)};
\node [below] at (2,-0.5) {\small odd($2$)};
\node [below] at (3,-0.5) {\small odd($3$)};
\node [below] at (4,-0.5) {\small odd($4$)};
\node [below] at (5,-0.5) {\small odd($5$)};

\node [left] at (-0.5,1) {\small even($1$)};
\node [left] at (-0.5,2) {\small even($2$)};
\node [left] at (-0.5,3) {\small even($3$)};
\node [left] at (-0.5,4) {\small even($4$)};
\node [left] at (-0.5,5) {\small even($5$)};
\end{tikzpicture}
\caption{Lattice of Pythagorean triples in $\gC$.}
\label{lattice1}
\end{figure}
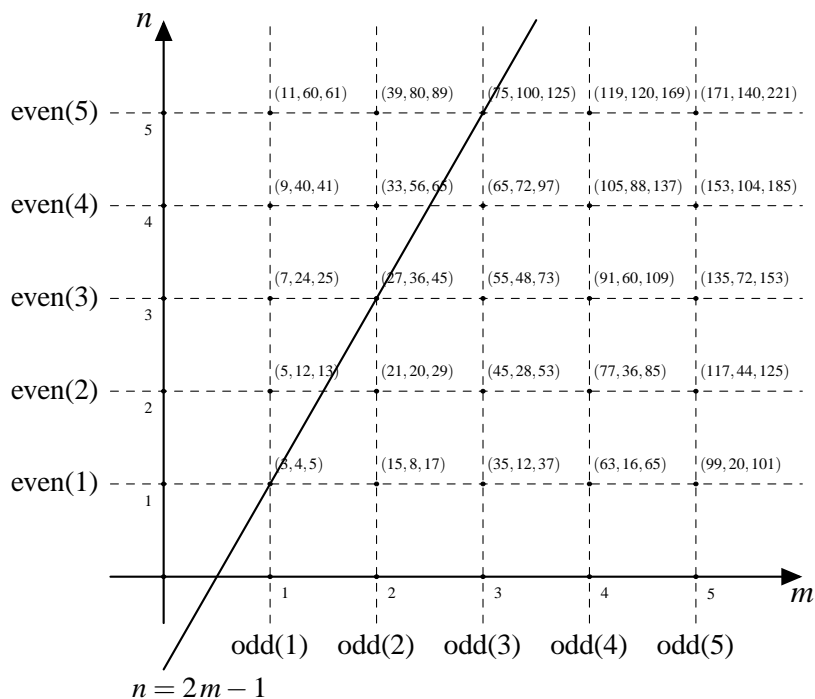
\end{center}
From Figure \ref{lattice1}, it is clear that $\gC$ also contains non-primitive triples. For example, all triples lying on the line $n=2\,m-1$ can be expressed as $p^2\,(3,4,5)$ where $p\in \mathbb{N}$ is odd. More in general, the following result holds.

\begin{theorem}
Let $(a,b,c) \in \gC$ and let $m,n\in \mathbb{N} \setminus \{0\}$ such that $(a,b,c) = \mathfrak{P}(m,n)$. Then, $(a,b,c)$ is primitive if and only if $n$ and $2\,m-1$ are relatively prime.
\end{theorem}
\proof
Let $\mu \defi 2\,m-1$. Writing (\ref{dioph1}-\ref{dioph3}) as a function of $n$ and $\mu$ yields
\bea
a \ns&\ns = \ns&\ns \mu\,(2\,n+\mu), \label{dioph1a} \\
b \ns&\ns = \ns&\ns 2\,n\,(n+\mu), \label{dioph2a} \\
c \ns&\ns = \ns&\ns 2\,n^2+\mu\,(2\,n+\mu).  \label{dioph3a}
\eea
Firstly, it is shown that if $n$ and $\mu$ have a common factor, $(a,b,c)$ is not primitive. Let us write $n=q\,\tilde{n}$ and $\mu=q\,\tilde{\mu}$, where $q,\tilde{n},\tilde{\mu} \in \mathbb{N} \setminus \{0\}$. Substituting these expressions into (\ref{dioph1a}-\ref{dioph3a}) yields
\beann
a \ns&\ns = \ns&\ns 2\,q^2\,\tilde{\mu}\,(2\,\tilde{n}+\tilde{\mu}), \\
b\ns&\ns = \ns&\ns  2\,q^2\,\tilde{n}\,(\tilde{n}+\tilde{\mu}), \\
c \ns&\ns = \ns&\ns q^2\,(2\,\tilde{n}^2+2\,\tilde{n}\,\tilde{\mu}+\tilde{\mu}^2),
\eeann
which show that $(a,b,c) \in \gC \setminus \gP_{\scriptscriptstyle 0}$. 

It is now shown that if $(a,b,c)$ is not primitive in $\gC$, then $n$ and $\mu$ are not relatively prime. If $(a,b,c)$ is not primitive, it can be written as $(a,b,c)=k\,(\alpha,\beta,\gamma)$ where $k \in \mathbb{N} \setminus \{0\}$ and $(\alpha,\beta,\gamma) \in \gP_{\scriptscriptstyle 0}$. Therefore, there exist $\hat{n},\hat{\mu} \in \mathbb{N} \setminus \{0\}$ such that $(\alpha,\beta,\gamma) =\mathfrak{P}(\hat{n},\hat{\mu})$, which, by taking (\ref{dioph1a}-\ref{dioph3a}) into account, gives
\beann
\mu\,(2\,n+\mu) \ns&\ns = \ns&\ns k\,\hat{\mu}\,(2\,\hat{n}+\hat{\mu}), \\
2\,n\,(n+\mu) \ns&\ns = \ns&\ns  2\,k\,\hat{n}\,(\hat{n}+\hat{\mu}), \\
2\,n^2+\mu\,(2\,n+\mu) \ns&\ns = \ns&\ns 2\,k\,\hat{n}^2+\hat{\mu}\,k\,(2\,\hat{n}+\hat{\mu}).
\eeann
The first and the third equations give $k=n^2 / \hat{n}^2$, which, once substituted into the second equation, yields
\[
\ell \defi \frac{n}{\hat{n}} = \frac{\mu}{\hat{\mu}}.
\] 
Notice that $\ell \in \mathbb{N}$ and $\ell^2 =k$. Thus, $\ell$ is a common factor of $n$ and $\mu$, which are therefore not relatively prime.
\endproof

\begin{remark}
{\em 
From the considerations in Remark \ref{remChat}, the set of Pythagorean triples in $\widehat{\gC}=\gE$ can be represented as the lattice in Figure \ref{fig3}.

\begin{center}
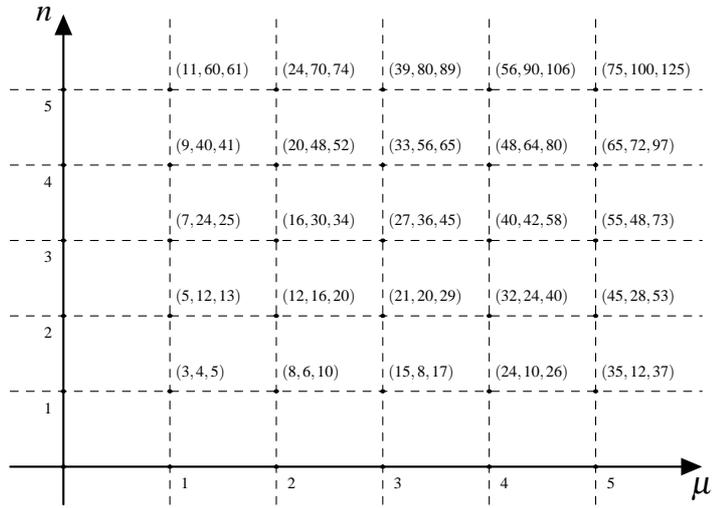
\begin{figure}[!h]
\hspace{4cm} \begin{tikzpicture}[line cap=round,line join=round,>=triangle 45,x=1.4cm,y=1.0cm]

\draw [->,thick] (-0.5,0) -- (6,0);
\draw [->,thick] (0,-0.5) -- (0,6);

\draw [fill] (0,0) circle [radius=0.02];
\draw [fill] (0,1) circle [radius=0.02];
\draw [fill] (0,2) circle [radius=0.02];
\draw [fill] (0,3) circle [radius=0.02];
\draw [fill] (0,4) circle [radius=0.02];
\draw [fill] (0,5) circle [radius=0.02];

\draw [fill] (1,0) circle [radius=0.02];
\draw [fill] (1,1) circle [radius=0.02];
\draw [fill] (1,2) circle [radius=0.02];
\draw [fill] (1,3) circle [radius=0.02];
\draw [fill] (1,4) circle [radius=0.02];
\draw [fill] (1,5) circle [radius=0.02];

\draw [fill] (2,0) circle [radius=0.02];
\draw [fill] (2,1) circle [radius=0.02];
\draw [fill] (2,2) circle [radius=0.02];
\draw [fill] (2,3) circle [radius=0.02];
\draw [fill] (2,4) circle [radius=0.02];
\draw [fill] (2,5) circle [radius=0.02];

\draw [fill] (3,0) circle [radius=0.02];
\draw [fill] (3,1) circle [radius=0.02];
\draw [fill] (3,2) circle [radius=0.02];
\draw [fill] (3,3) circle [radius=0.02];
\draw [fill] (3,4) circle [radius=0.02];
\draw [fill] (3,5) circle [radius=0.02];

\draw [fill] (4,0) circle [radius=0.02];
\draw [fill] (4,1) circle [radius=0.02];
\draw [fill] (4,2) circle [radius=0.02];
\draw [fill] (4,3) circle [radius=0.02];
\draw [fill] (4,4) circle [radius=0.02];
\draw [fill] (4,5) circle [radius=0.02];

\draw [fill] (5,0) circle [radius=0.02];
\draw [fill] (5,1) circle [radius=0.02];
\draw [fill] (5,2) circle [radius=0.02];
\draw [fill] (5,3) circle [radius=0.02];
\draw [fill] (5,4) circle [radius=0.02];
\draw [fill] (5,5) circle [radius=0.02];

%
\draw [dashed,thin] (-0.5,1) -- (6,1);
\draw [dashed,thin] (-0.5,2) -- (6,2);
\draw [dashed,thin] (-0.5,3) -- (6,3);
\draw [dashed,thin] (-0.5,4) -- (6,4);
\draw [dashed,thin] (-0.5,5) -- (6,5);

\draw [dashed,thin] (1,-0.5) -- (1,6);
\draw [dashed,thin] (2,-0.5) -- (2,6);
\draw [dashed,thin] (3,-0.5) -- (3,6);
\draw [dashed,thin] (4,-0.5) -- (4,6);
\draw [dashed,thin] (5,-0.5) -- (5,6);
\node [above right] at (1,1) {\tiny $\!\!\!(3,4,5)$};
\node [above right] at (2,1) {\tiny $\!\!\!(8,6,10)$};
\node [above right] at (3,1) {\tiny $\!\!\!(15,8,17)$};
\node [above right] at (4,1) {\tiny $\!\!\!(24,10,26)$};
\node [above right] at (5,1) {\tiny $\!\!\!(35,12,37)$};

\node [above right] at (1,2) {\tiny $\!\!\!(5,12,13)$};
\node [above right] at (2,2) {\tiny $\!\!\!(12,16,20)$};
\node [above right] at (3,2) {\tiny $\!\!\!(21,20,29)$};
\node [above right] at (4,2) {\tiny $\!\!\!(32,24,40)$};
\node [above right] at (5,2) {\tiny $\!\!\!(45,28,53)$};

\node [above right] at (1,3) {\tiny $\!\!\!(7,24,25)$};
\node [above right] at (2,3) {\tiny $\!\!\!(16,30,34)$};
\node [above right] at (3,3) {\tiny $\!\!\!(27,36,45)$};
\node [above right] at (4,3) {\tiny $\!\!\!(40,42,58)$};
\node [above right] at (5,3) {\tiny $\!\!\!(55,48,73)$};

\node [above right] at (1,4) {\tiny $\!\!\!(9,40,41)$};
\node [above right] at (2,4) {\tiny $\!\!\!(20,48,52)$};
\node [above right] at (3,4) {\tiny $\!\!\!(33,56,65)$};
\node [above right] at (4,4) {\tiny $\!\!\!(48,64,80)$};
\node [above right] at (5,4) {\tiny $\!\!\!(65,72,97)$};

\node [above right] at (1,5) {\tiny $\!\!\!(11,60,61)$};
\node [above right] at (2,5) {\tiny $\!\!\!(24,70,74)$};
\node [above right] at (3,5) {\tiny $\!\!\!(39,80,89)$};
\node [above right] at (4,5) {\tiny $\!\!\!(56,90,106)$};
\node [above right] at (5,5) {\tiny $\!\!\!(75,100,125)$};

\node [below] at (6,0) { $\mu$};
\node [left] at (0,6) { $n$};

\node [below right] at (1,0) {\tiny $1$};
\node [below right] at (2,0) {\tiny $2$};
\node [below right] at (3,0) {\tiny $3$};
\node [below right] at (4,0) {\tiny $4$};
\node [below right] at (5,0) {\tiny $5$};

\node [below left] at (0,1) {\tiny $1$};
\node [below left] at (0,2) {\tiny $2$};
\node [below left] at (0,3) {\tiny $3$};
\node [below left] at (0,4) {\tiny $4$};
\node [below left] at (0,5) {\tiny $5$};
\end{tikzpicture}
\caption{Lattice of Pythagorean triples in $\widehat{\gC}=\gE$.}
\label{fig3}
\end{figure}
\end{center}
}
\end{remark}

\begin{remark}
{\em
Particularly important special cases of the parameterisation given in Theorem \ref{cor} are:
\begin{enumerate}
\item the so-called {\em Pythagorean family} of odd triples $(2\,n+1,2\,n^2+2\,n,2\,n^2+2\,n+1)$, which can be obtained as $\mathfrak{P}(1,n)=\text{odd}(1)=\{(3,4,5),(5,12,13),(7,24,25)\ldots\}$, i.e.
\[
\mathfrak{P}(1,n)=\left\{ (a,b,c)\in \mathbb{N}^3\; : \;\;\bmat{c} a \\ b \\ c \emat=\bmat{ccc} 0 & 1 & 1 \\ 1 & 1 & 0 \\ 1 & 1 & 1 \emat \bmat{c} 2\,n^2 \\ 2\,n \\ 1 \emat \text{ for some } n \in \mathbb{N} \setminus \{0\}\right\},
\]
which follows also from (\ref{matrixform}) by setting $m=1$;
\item the so-called {\em Platonic family} of even triples $(4\,m^2-1,4\,m,4\,m^2+1)$, which can be obtained as $\mathfrak{P}(m,1)=\text{even}(1)=\{(3,4,5),(15,8,17),(35,12,37)\ldots\}$, i.e.
\[
\mathfrak{P}(m,1)=\left\{ (a,b,c)\in \mathbb{N}^3\; : \;\;\bmat{c} a \\ b \\ c \emat=\bmat{ccc} 0 & 1 & 1 \\ 1 & 1 & 0 \\ 1 & 1 & 1 \emat \bmat{c} 2 \\ 2\,(2\,m-1) \\ (2\,m-1)^2 \emat \text{ for some } m \in \mathbb{N} \setminus \{0\}\right\}
\]
which follows from (\ref{matrixform}) with $n=1$.
\end{enumerate}
The only Pythagorean triple that simultaneously belongs to the set of Pythagorean family of triples and to the Platonic family of triples is obviously the Pythagorean triple $(3,4,5)$. There holds $\mathfrak{P}(1,1)=(3,4,5)$.
\

}
\end{remark}

\section*{Concluding remarks}
{
In this paper  an alternative parameterisation of Pythagorean triples has been introduced, which hinges on the series of odd and even triples. This parameterisation generates all the primitive Pythagorean triples, but, with respect to the set of triples that can be generated using Euclid's formula, does not include the non-primitive Pythagorean triples $(a,b,c)$ where $a$, $b$ and $c$ are all even. 
}

\end{document}